# LARGE DEVIATIONS FOR TEMPLATE MATCHING BETWEEN POINT PROCESSES

By Zhiyi Chi

*University of Chicago*

We study the asymptotics related to the following matching criteria for two independent realizations of point processes $X \sim \mathbf{X}$ and $Y \sim \mathbf{Y}$. Given $l > 0$, $X \cap [0, l)$ serves as a template. For each $t > 0$, the matching score between the template and $Y \cap [t, t+l)$ is a weighted sum of the Euclidean distances from $y - t$ to the template over all $y \in Y \cap [t, t+l)$. The template matching criteria are used in neuroscience to detect neural activity with certain patterns. We first consider $W_l(\theta)$, the waiting time until the matching score is above a given threshold $\theta$. We show that whether the score is scalar- or vector-valued, $(1/l) \log W_l(\theta)$ converges almost surely to a constant whose explicit form is available, when $\mathbf{X}$ is a stationary ergodic process and $\mathbf{Y}$ is a homogeneous Poisson point process. Second, as $l \to \infty$, a strong approximation for $-\log[\Pr\{W_l(\theta) = 0\}]$ by its rate function is established, and in the case where $\mathbf{X}$ is sufficiently mixing, the rates, after being centered and normalized by $\sqrt{l}$, satisfy a central limit theorem and almost sure invariance principle. The explicit form of the variance of the normal distribution is given for the case where $\mathbf{X}$ is a homogeneous Poisson process as well.

**1. Introduction.** In neuroscience, it is well accepted that neurons are the basic units of information processing. By complex biochemical mechanisms governing the ion flows through its membrane, a neuron generates very narrow and highly peaked electric potentials, or "spikes," in its soma (main body) [6]. These spikes can propagate along the neuron's axons, which are cables that extend over relatively long distance to reach the other cells. The spikes can then influence the activities of those cells. The temporal pattern in which a neuron generates spikes dynamically depends on its inputs, which are either stimuli from the environment or biochemicals induced by the spikes from the other neurons. In this way, information is processed through the neural network. Because spikes are very narrow and peaked, point processes are the most commonly used models for neuronal activity, with points representing the temporal locations of spikes.

For many studies in neuroscience, it is necessary to detect segments of neuronal activity that exhibit certain patterns [1, 10, 11]. Recently, in a study on the activity of brain during sleep, a template matching algorithm was developed which uses linear filtering to quickly detect such segments (cf. [3]). The algorithm is template based. Suppose $S = \{x_1, \ldots, x_n\}$ is a nonempty sequence of spikes generated by a neuron under some specific condition between time 0 and $l$. This sequence is used as a template. Given a data sequence of spikes $Y = \{y_1, y_2, \ldots\}$ generated by the same neuron but at a different time, the goal is to find segments in $Y$ that have a temporal pattern similar to $S$. To do this, for each time point $t$, collect all $y$'s between $t$ and $t+l$ and shift them back to the origin. If the temporal distances between the shifted $y$'s and $S$ are small on average, then it indicates that the temporal pattern of the activity recorded in $Y$ between $t$ and $t+l$ is similar to that of $S$. Therefore one can use the following matching score

$$M(t) = \frac{1}{l} \sum_{y \text{ between } t \text{ and } t+l} f(d(y-t, S))$$

to measure the overall distance, where $f(x)$ is a function of $x \geq 0$ that is nonincreasing, and $d$ is the Euclidean distance such that for any $y \in \mathbb{R}$ and $S \subset \mathbb{R}$, $d(y, S) = \inf\{|y - s| : s \in A\}$. Let $\theta$ be a threshold value fixed beforehand. If $M(t) \geq \theta$, then output $t$ as a location of matching segment, or "target." To improve accuracy, the detection was modified to involve multiple matching criteria so that both $f$ and $\theta$ are vector-valued. Then $t$ is a target location only if $M(t) \geq \theta$ (cf. [3]), where, for $u = (u_1, \ldots, u_n)$ and $v = (v_1, \ldots, v_n)$, "$u \geq v$" denotes "$u_j \geq v_j$ for all $j$." For later use, let "$u > v$" denote "$u \geq v$ and $u \neq v$."

In the above studies, it is necessary to evaluate how difficult it is to get false targets if a data sequence is noise. A useful criterion for this is the waiting time until the matching score is larger than or equal to $\theta$. Presumably, when the template is longer, that is, $l$ is larger, it would be more difficult to find false targets. But how much more difficult? In this article, we study the asymptotics of the waiting time under certain assumptions on the point processes underlying the template and the data.

To fix notation, realizations of a point process on $\mathbb{R}$ will be regarded as point sequences. For $a < b$ and $S \subset \mathbb{R}$, denote

$$S_a^b = S \cap [a, b), \qquad S - a = \{t - a : a \in S\}.$$

We will think of the template $S$ as an initial segment of an infinite sequence $X$ of points on $\mathbb{R}$. That is, $S = X_0^l$ for some $l > 0$. Given $f =$ $(f_1, \ldots, f_n) : \{0\} \cup \mathbb{R} \to \mathbb{R}^n$, if $Y$ is another sequence of points, then for each $t > 0$, define

$$\rho_l(X_0^l, Y_t^{t+l}) = \begin{cases} \dfrac{1}{l} \displaystyle\sum_{y \in Y_t^{t+l}} f(d(y - t, X_0^l)), & \text{if } X_0^l \neq \varnothing, \\ (-\infty, \ldots, -\infty), & \text{otherwise.} \end{cases}$$

In practice, it is reasonable to require that $f_k(x)$, $k = 1, \ldots, n$, be non-increasing functions in $x \geq 0$. However, to get the asymptotics of $W$, this requirement can be dropped. Given a threshold $\theta \in \mathbb{R}^n$, the waiting time until the first false target is detected is

$$W_l(\theta, X, Y) = \inf\{t \geq 0 : \rho_l(X_0^l, Y_t^{t+l}) \geq \theta\}.$$

To study the asymptotics of $W_l$ as $l$ increases, assume $X$ and $Y$ are random realizations of two point processes $\mathbf{X}$ and $\mathbf{Y}$ on $\mathbb{R}$, respectively. One would think of stationary Poisson point processes as signals that contain the least amount of information. In other words, they are plainly noise. We will mainly focus on the case where $\mathbf{Y}$ is Poisson.

The asymptotics of waiting times for pattern detection using random templates have been studied for the case where $\mathbf{X} = \{X_n, n \geq 1\}$ and $\mathbf{Y} = \{Y_n, n \geq 1\}$ are integer indexed processes (cf. [2, 7, 13, 14] and references therein). In these works, the matching score is defined for $(X_1, \ldots, X_n)$ and $(Y_1, \ldots, Y_n)$ as the average of $\rho(X_j, Y_j)$ for some function $\rho$. Whereas the temporal relations between points are essential in the asymptotics considered here, it is apparent such relations are not relevant in the above results.

When $f$ is scalar-valued function $f$, the first main result is:

THEOREM 1. *Suppose that $\mathbf{X}$ and $\mathbf{Y}$ are point processes on $\mathbb{R}$ that are independent of each other and $f$ is a bounded scalar function. Assume:*

1. $\mathbf{X}$ *is a stationary and ergodic point process with mean density*

$$\overline{N} = EN_{\mathbf{X}}[0, 1) \in (0, \infty),$$

*where $N_{\mathbf{X}}(\cdot)$ is the random counting measure associated with $\mathbf{X}$ (cf. [5]).*
2. $\Pr\{d(0, \mathbf{X})$ *is a continuity point of* $f\} = 1$.
3. $\Pr\{f(d(0, \mathbf{X})) > 0\} > 0$.
4. $\mathbf{Y}$ *is a Poisson point process with density* $\lambda \in (0, \infty)$.

*Define*

(1.1) $$\phi := \lambda E[f(d(0, \mathbf{X}))],$$

(1.2) $$\Lambda(t) := \lambda E[e^{tf(d(0, \mathbf{X}))} - 1].$$

*Then, given $\theta > \phi$,*

(1.3) $$\lim_{l \to \infty} \frac{1}{l} \log W_l(\theta, \mathbf{X}, \mathbf{Y}) = \sup_{t \geq 0}\{\theta t - \Lambda(t)\} \quad w.p.1.$$

Theorem 1 can be generalized to the case where the signal is a compound Poisson process. Such a process can be characterized as a pair $\widetilde{\mathbf{Y}} = (\mathbf{Y}, \{Q(y), y \in \mathbb{R}\})$, where $\mathbf{Y}$ is a common Poisson point process with density $\lambda$ and $Q(y)$ i.i.d. $\sim Q \in \mathbb{N}$ are random variables independent of $\mathbf{X}$ and $\mathbf{Y}$. For $Y \sim \mathbf{Y}$, each $y \in Y$ is interpreted as a location where there is at least one point, and $Q(y)$ is the number of points at $y$. Then for $\widetilde{Y} \sim \widetilde{\mathbf{Y}}$, the matching score between $X_0^l$ and the segment of $\widetilde{Y}$ in $[t, t+l)$, denoted by $\widetilde{Y}_t^{t+l}$, is

$$\rho_l(X_0^l, \widetilde{Y}_t^{t+l}) = \frac{1}{l} \sum_{y \in Y_t^{t+l}} Q(y) f(d(y-t, X_0^l)).$$

PROPOSITION 1. *Suppose all the assumptions in Theorem 1 are satisfied. In addition, suppose $G(t) := E[e^{tQ}] < \infty$ for all $t > 0$. Then, given $\theta > \phi := \lambda E[f(d(0, \mathbf{X}))] E[Q]$,*

$$\lim_{l \to \infty} \frac{1}{l} \log W_l(\theta, \mathbf{X}, \widetilde{\mathbf{Y}}) = \sup_{t \geq 0} \{\theta t - \tilde{\Lambda}(t)\} \qquad w.p.1,$$

*where $\tilde{\Lambda}(t) = \lambda E[G(tf(d(0, \mathbf{X}))) - 1]$.*

The asymptotic in Theorem 1 can also be proved when $n = \dim f > 1$. Because the monotonicity property of $\mathbb{R}$ used in the proof of Theorem 1 is lost in this case, some changes in the assumptions are needed.

THEOREM 2. *Assume $\mathbf{X}$, $\mathbf{Y}$ and $f$ satisfy all but condition 3 in Theorem 1. Instead, assume:*

3′. *For any $v \neq 0$, $\Pr\{\langle v, f(d(0, \mathbf{X}))\rangle > 0\} > 0$.*

*Define $\Lambda(t) = \lambda E[e^{\langle t, f(d(0, \mathbf{X}))\rangle} - 1]$ and $\phi$ as in (1.1). Then for any $\theta > \phi$,*

$$(1.4) \qquad \lim_{l \to \infty} \frac{1}{l} \log W_l(\theta, \mathbf{X}, \mathbf{Y}) = \inf_{z \geq \theta} \Lambda^*(z) \qquad w.p.1,$$

*where*

$$\Lambda^*(z) = \sup_{t \in \mathbb{R}^n} \{\langle z, t\rangle - \Lambda(t)\}$$

*is bounded and continuous.*

The proofs for Theorems 1 and 2 rely on the conditional large deviations principle (LDP) of a family of random variables, because $X \sim \mathbf{X}$ is a fixed realization (cf. [2, 4, 7, 8]). These random variables have close relationship to $\rho_l(X_0^l, \mathbf{Y}_0^l)$. We next consider the asymptotics of the latter and restrict our focus to the case where $f$ is scalar-valued. First, the following approximation for the conditional LDP for $\rho_l(X_0^l, \mathbf{Y}_0^l)$ holds.

THEOREM 3. *Under the same assumption as in Theorem 1, for any set of points $S \subset \mathbb{R}$, let*

$$(1.5) \quad \Lambda_{S,l}(t) = \begin{cases} \dfrac{1}{l} \log E\left[\exp\left\{t \sum_{y \in \mathbf{Y}_0^l} f(d(y,S))\right\}\right] = \dfrac{\lambda}{l} \int_0^l [e^{tf(d(y,S))} - 1] \, dy, \\ \qquad \qquad \qquad \qquad \qquad \qquad \qquad \qquad \qquad \qquad \text{if } S \neq \varnothing, \\ 0, \qquad \qquad \qquad \qquad \qquad \qquad \qquad \qquad \qquad \qquad \text{otherwise,} \end{cases}$$

$$(1.6) \quad \Lambda^*_{S,l}(\theta) = \sup_{t \in \mathbb{R}}[\theta t - \Lambda_{S,l}(t)].$$

*Then, given $\theta > \phi$, almost surely, for $X \sim \mathbf{X}$,*

$$(1.7) \qquad -\log \Pr\{\rho_l(X_0^l, \mathbf{Y}_0^l) \geq \theta\} - l\Lambda^*_{X_0^l, l}(\theta) = o(\sqrt{l}).$$

*Note that* $\Pr\{\rho_l(X_0^l, \mathbf{Y}_0^l) \geq \theta\} = \Pr\{W_l(\theta) = 0\}$.

REMARK. Despite the higher-order approximation in Theorem 3, the difference between the aforementioned random variables and $\rho_l(X_0^l, \mathbf{Y}_0^l)$ does not allow the approximation to be applied to the proof of Theorem 1 and it is not clear to me how to derive a similar higher-order approximation to $W_l$.

Finally, under suitable conditions, $-\log \Pr\{\rho_l(X_0^l, \mathbf{Y}_0^l) \geq \theta\}$ after being centered and normalized is asymptotically normal, as the following result combined with (1.7) shows.

THEOREM 4. *Assume $\mathbf{X}$, $\mathbf{Y}$ and $f$ satisfy all but condition 2 in Theorem 1. Instead, assume $f \neq 0$ is continuous. Given $\theta > \phi$, let $t_0$ be the (unique) point with $\Lambda^*(\theta) = \theta t_0 - \Lambda(t_0)$. If $\mathbf{X}$ is a Poisson point process with density $\rho$, then almost surely, for $X \sim \mathbf{X}$,*

$$(1.8) \qquad l\{\Lambda^*_{X_0^l, l}(\theta) - [\theta t_0 - \Lambda_{X_0^l, l}(t_0)]\} = o(\sqrt{l}), \qquad l \to \infty,$$

$$(1.9) \qquad \sqrt{l}(\theta t_0 - \Lambda_{X_0^l, l}(t_0) - \Lambda^*(\theta)) \xrightarrow{D} N(0, 4\rho\sigma^2),$$

*with*

$$(1.10) \quad \sigma^2 = \operatorname{Var}\left\{G\left(\dfrac{U}{2\rho}\right) - UE\left[g\left(\dfrac{U}{2\rho}\right)\dfrac{U}{2\rho}\right]\right\} + \left\{E\left[G\left(\dfrac{U}{2\rho}\right) - g\left(\dfrac{U}{2\rho}\right)\dfrac{U}{2\rho}\right]\right\}^2,$$

*where $U \sim \operatorname{Exp}(1)$, $g(x) = e^{t_0 f(x)}$ and $G(x) = \int_0^x g$.*

REMARK. Following the proof of Theorem 4, it can be shown that, instead of assuming $\mathbf{X}$ to be a Poisson process, if $\int_0^\infty \psi(t) \, dt < \infty$ and either $f$ has bounded support or $E\tau^2 < \infty$, where $\psi(t) = \sup\{|P(A \cap B) - P(A)P(B)| : A \in \sigma(\mathbf{X}_{-\infty}^0), B \in \sigma(\mathbf{X}_t^\infty)\}$ and $\tau = \min(\mathbf{X}_0^\infty)$, then (1.8) and

the asymptotic normality of $\sqrt{l}(\theta t_0 - \Lambda_{X_0^l,l}(t_0) - \Lambda^*(\theta))$ still hold. Indeed, under the assumptions, the left-hand side of (1.9) is $\sqrt{n}(\Lambda_{X,n}(t_0) - \Lambda(t_0)) + o(1)$, w.p.1, with $n = \lfloor l \rfloor$, and the random variables $Z_n = \lambda \int_n^{n+1} [e^{t_0 f(d(y, \mathbf{X}))} - 1] \, dy$ satisfy the mixing condition in [12], Theorem 1, yielding the asymptotic normality. However, in general, the explicit form of the limit distribution is not readily obtained.

The rest of the article is organized as follows. In Sections 2 and 3, Theorem 1 is proved. In Section 4, Theorem 2 is proved. In Section 5, Theorem 3 is proved. Finally, in Section 6, Theorem 4 is proved.

**2. Waiting times for scalar-valued matching scores.** In this section, suppose $\mathbf{X}$ and $\mathbf{Y}$ satisfy the conditions in Theorem 1. For any function $g$, denote $g^+ = \max(g, 0)$ and $g^- = \max(-g, 0)$, and for $\varepsilon > 0$,

$$g_\varepsilon(x) = \sup_{|t-x| \leq \varepsilon} g(t).$$

For integer $n > 1$ and $X, Y \subset \mathbb{R}$ with $Y$ discrete, define

$$A_n(X, Y) = \frac{1}{n} \sum_{y \in Y_0^{n-1}} \inf_{n-1 \leq l \leq n} f^+(d(y, X_0^l))$$

$$- \frac{1}{n-1} \sum_{y \in Y_0^n} \sup_{n-1 \leq l \leq n} f^-(d(y, X_0^l)),$$

$$B_{n,\varepsilon}(X, Y) = \frac{1}{n-1} \sum_{y \in Y_0^{n+\varepsilon}} \sup_{n-1 \leq l \leq n} f_\varepsilon^+(d(y, X_0^l))$$

$$- \frac{1}{n} \sum_{y \in Y_\varepsilon^{n-1}} \inf_{n-1 \leq l \leq n} f_\varepsilon^-(d(y, X_0^l)).$$

Since

$$(2.1) \quad f_\varepsilon^+(x) = \sup_{|t-x| \leq \varepsilon} f^+(t) \geq f^+(x), \qquad f_\varepsilon^-(x) = \inf_{|t-x| \leq \varepsilon} f^-(t) \leq f^-(x),$$

it is seen that $B_{n,\varepsilon}(X, Y) \geq A_n(X, Y)$. The following lemmas are needed for the proof of Theorem 1.

LEMMA 1. *Given $\theta \in \mathbb{R}$, almost surely, for $X \sim \mathbf{X}$, as $n \to \infty$, eventually there are*

$$\alpha_{n,X,\theta} := \Pr\{A_n(X, \mathbf{Y}) \geq \theta\} > 0, \qquad \beta_{n,\varepsilon,X,\theta} := \Pr\{B_{n,\varepsilon}(X, \mathbf{Y}) \geq \theta\} > 0.$$

Because of Lemma 1, the logarithms in the results below are well defined almost surely.

LEMMA 2 (Upper bounds for $W_l$). *Let $\theta$ be an arbitrary number. Then*

$$\text{(2.2)} \quad \Pr\left\{\limsup_{l\to\infty} \frac{1}{l} \log[\alpha_{\lceil l\rceil, \mathbf{X}, \theta} \times W_l(\theta, \mathbf{X}, \mathbf{Y})] \leq 0\right\} = 1.$$

LEMMA 3 (Lower bounds for $W_l$). *Let $\theta$ be an arbitrary number. Then*

$$\text{(2.3)} \quad \Pr\left\{\liminf_{l\to\infty} \frac{1}{l} \log[\beta_{\lceil l\rceil, \varepsilon, \mathbf{X}, \theta} \times \max(W_l(\theta, \mathbf{X}, \mathbf{Y}), 1)] \geq 0\right\} = 1.$$

LEMMA 4 (LDP). *Almost surely, for $X \sim \mathbf{X}$, the conditional laws of $A_n(X, \mathbf{Y})$, $n \geq 2$, satisfy the LDP with a good rate function*

$$\text{(2.4)} \quad \Lambda^*(\theta) = \sup_{t\in\mathbb{R}}\{\theta t - \Lambda(t)\},$$

*and the conditional laws of $B_{n,\varepsilon}(X, \mathbf{Y})$, $n \geq 2$, satisfy the LDP with a good rate function*

$$\Lambda_\varepsilon^*(\theta) = \sup_{t\in\mathbb{R}}\{\theta t - \lambda E[e^{tf_\varepsilon(d(0,\mathbf{X}))} - 1]\}.$$

Assume for now that the above lemmas hold. For $\theta > \phi$, by Lemmas 2 and 4, almost surely, for $X \sim \mathbf{X}, Y \sim \mathbf{Y}$,

$$\text{(2.5)} \quad \limsup_{l\to\infty} \frac{1}{l} \log W_l(\theta, X, Y) \leq \inf_{z>\theta} \Lambda^*(z).$$

It is known that $\Lambda$ is strictly convex (e.g., [9]). Because $f$ is bounded, $\Lambda$ is smooth everywhere with $\Lambda'(0) = \phi$. By condition 3 of Theorem 1, $\Lambda(t) \to \infty$ exponentially as $t \to \infty$. These imply that for any $z > \phi$, $\Lambda^*(z) > 0$ is finite and achieved on $(0, \infty)$, and $\Lambda^*$ is a continuous strictly increasing convex function on $(\phi, \infty)$. Then by (2.5), it is seen that

$$\text{(2.6)} \quad \limsup_{l\to\infty} \frac{1}{l} \log W_l(\theta, X, Y) \leq \Lambda^*(\theta),$$

and to complete the proof of (1.3), it remains to show

$$\text{(2.7)} \quad \liminf_{l\to\infty} \frac{1}{l} \log W_l(\theta, X, Y) \geq \Lambda^*(\theta).$$

By Lemmas 3 and 4, for any $\varepsilon > 0$,

$$\liminf_{l\to\infty} \frac{1}{l} \log \max(W_l(\theta, X, Y), 1) \geq \inf_{z\geq\theta} \Lambda_\varepsilon^*(z).$$

Similar to the above argument, it is seen that almost surely, for $X \sim \mathbf{X}$, $Y \sim \mathbf{Y}$,

$$\text{(2.8)} \quad \liminf_{l\to\infty} \frac{1}{l} \log \max(W_l(\theta, X, Y), 1) \geq \Lambda_\varepsilon^*(\theta) = \sup_{t\geq 0}\{\theta t - \Lambda_\varepsilon(t)\},$$

where $\Lambda_\varepsilon(t) = \lambda E[e^{tf_\varepsilon(d(0,\mathbf{X}))} - 1]$. Let $t^*$ be the unique point where $\Lambda^*(\theta) = \theta t^* - \Lambda(t^*)$. Then $\Lambda^*_\varepsilon(\theta) \geq \theta t^* - \Lambda_\varepsilon(t^*)$. By condition 2 of Theorem 1 and dominated convergence, $\Lambda_\varepsilon(t^*) \to \Lambda(t^*)$, leading to $\liminf_{\varepsilon \to 0} \Lambda^*_\varepsilon(\theta) \geq \Lambda^*(\theta) > 0$. So by (2.8)

$$\liminf_{l \to \infty} \frac{1}{l} \log \max(W_l(\theta, X, Y), 1) \geq \Lambda^*(\theta) > 0.$$

The lower bound also implies $W_l(\theta, X, Y) \to \infty$. These combined with (2.8) prove (2.7). □

## 3. Proofs of lemmas.

PROPOSITION 2. *For $\mathbf{X}$ satisfying condition 1 of Theorem 1,*

$$\Pr\left\{\lim_{l \to \infty} \frac{l - \sup\{x : x \in \mathbf{X}_0^l\}}{l} = 0\right\} = 1,$$

*where, for $\mathbf{X}_0^l = \varnothing$, $\sup\{x : x \in \mathbf{X}_0^l\}$ is defined to be $-\infty$.*

PROOF. Because $\mathbf{X}$ is stationary and ergodic, almost surely, for a realization $X$ of $\mathbf{X}$, as $l \to \infty$, $N_X[0, l]/l \to \overline{N} > 0$, implying that for any $\varepsilon \in (0,1)$, $N_X[(1-\varepsilon)l, l) \to \infty$. Now

$$\frac{l - \sup\{X_0^l\}}{l} \geq \varepsilon \implies N_X((1-\varepsilon)l, l) = 0,$$

leading to

$$\Pr\left\{\limsup_{l \to \infty} \frac{l - \sup\{x : x \in \mathbf{X}_0^l\}}{l} \geq \varepsilon\right\} = 0,$$

which completes the proof. □

PROOF OF LEMMA 1. Because $B_{n,\varepsilon} \geq A_n$, it is enough to show that almost surely, for $X \sim \mathbf{X}$, $\alpha_{n,X,\theta} := \Pr\{A_n(X, \mathbf{Y}) \geq \theta\} > 0$ eventually, as $n \to \infty$. Let $X$ be a realization of $\mathbf{X}$ and $s_n = \min(X_0^{n-1})$, $\tau_n = \max(X_0^{n-1})$, for $n \geq 2$. It is easy to see $s_n/n \to 0$ w.p.1. By Proposition 2, almost surely, $\tau_n$ is well defined for all large $n$ and $(n - \tau_n)/n \to 0$. Note that for $y \in Y_{s_n}^{\tau_n}$, $d(y, X_0^{n-1}) = d(y, X)$. By the ergodicity of $\mathbf{X}$ and condition 3 of Theorem 1, almost surely,

$$\lim_{n \to \infty} E\left[\frac{1}{n} \sum_{y \in \mathbf{Y}_{s_n}^{\tau_n}} \mathbf{1}_{\{f(d(y, X_0^{n-1})) > 0\}}\right] = \lim_{n \to \infty} \frac{\lambda}{n} \int_{s_n}^{\tau_n} \mathbf{1}_{\{f(d(y,X)) > 0\}} \, dy$$

$$= \lim_{n \to \infty} \frac{\lambda}{n} \int_0^n \mathbf{1}_{\{f(d(0, X-y)) > 0\}} \, dy$$

$$= \lambda \Pr\{f(d(0, \mathbf{X})) > 0\} > 0.$$

Then it is seen that for $n$ large enough, there is $\eta_n > 0$ such that

$$\Pr\left\{\frac{1}{n}\sum_{y \in \mathbf{Y}_0^{\tau_n}} f^+(d(y, X_0^{n-1})) > \eta_n\right\} > 0.$$

Define

$$C_n = \left\{Y : Y_0^n = Y_{s_n}^{\tau_n}, \frac{1}{n}\sum_{y \in Y_0^n} f(d(y, X_0^{n-1})) > \eta_n \right.$$

$$\left. \text{and } \forall y \in Y_0^n, f(d(y, X_0^{n-1})) > 0\right\}.$$

By the property of Poisson processes, it is not hard to see that $\Pr\{\mathbf{Y} \in C_n\} > 0$. Fix $N \in \mathbb{N}$ with $N > \theta/\eta_n$. Let $D_n$ consist of all $Y$ with $Y_0^n$ being the union of $Z_1 \cap [0, n), \ldots, Z_N \cap [0, n)$ for some $Z_1, \ldots, Z_N \in C_n$ with $Z_i \cap Z_j \cap [0, n) = \varnothing$, $i \neq j$. Then $\Pr\{\mathbf{Y} \in D_n\} > (\Pr\{\mathbf{Y} \in C_n\})^N > 0$ and for any $Y \in D_n$, $A_n(X, Y) \geq N\eta_n \geq \theta$. $\square$

PROOF OF LEMMA 2. Let $\{K_n\}$ be a sequence of positive numbers to be determined later. Fix $n \geq 2$. Let $X$ be a realization of $\mathbf{X}$ with $\alpha_{n,X,\theta} > 0$ and $Y$ a realization of $\mathbf{Y}$. If there is $l \in (n-1, n]$, such that $W_l(\theta, X, Y) > K_n$, then for all $t \in [0, K_n]$,

$$\frac{1}{l}\sum_{y \in Y_t^{t+l}} f(d(y-t, X_0^l)) < \theta$$

$$\implies \frac{1}{l}\sum_{y \in Y_t^{t+l}} f^+(d(y-t, X_0^l))$$

$$< \theta + \frac{1}{l}\sum_{y \in Y_t^{t+l}} f^-(d(y-t, X_0^l))$$

$$\underset{(a)}{\implies} \frac{1}{n}\sum_{y \in Y_t^{t+n-1}} \inf_{n-1 \leq l \leq n} f^+(d(y-t, X_0^l))$$

$$< \theta + \frac{1}{n-1}\sum_{y \in Y_t^{t+n}} \sup_{n-1 \leq l \leq n} f^-(d(y-t, X_0^l))$$

$$\implies A_n(X, Y - t) < \theta,$$

with $(a)$ due to $f^+, f^- \geq 0$. In particular, $W_l(\theta, X, Y) > K_n$ implies $A_n(X, Y - kn) < \theta$ for $k = 0, \ldots, \lfloor K_n/n \rfloor$. Because $A_n(X, Y - kn)$ only depends on $X$ and $Y_{kn}^{(k+1)n}$, by the fact that $\mathbf{Y}_{kn}^{(k+1)n}$ are i.i.d.,

$$\Pr\{\exists l \in (n-1, n] \text{ s.t. } W_l(\theta, X, \mathbf{Y}) > K_n\}$$

$$\leq \Pr\left\{\bigcap_{k=0}^{\lfloor K_n/n \rfloor} \{A_n(X, \mathbf{Y} - kn) < \theta\}\right\}$$

$$= \prod_{k=0}^{\lfloor K_n/n \rfloor} \Pr\{A_n(X, \mathbf{Y} - kn) < \theta\} \leq (1 - \alpha_{n,X,\theta})^{K_n/n} \leq e^{-\alpha_{n,X,\theta} K_n/n}.$$

Choose $K_n = c(n)n/\alpha_{n,X,\theta}$, with $\sum e^{-c(n)} < \infty$ and $\frac{1}{n} \log c(n) \to 0$. Then

$$\Pr\left\{\exists l \in (n-1, n] \text{ s.t. } \frac{1}{l} \log[\alpha_{n,X,\theta} \times W_l(\theta, X, \mathbf{Y})] > \frac{1}{l} \log[c(n)n]\right\} \leq e^{-c(n)}.$$

Because the above bound is uniform over $X$ with $\alpha_{n,X,\theta} > 0$ and summable, by the Borel–Cantelli lemma and Lemma 1, (2.2) is therefore proved. □

PROOF OF LEMMA 3. Fix $n \geq 2$, $\varepsilon \in (0,1)$ and $L > 0$. Let $X$ be a realization of $\mathbf{X}$ with $\beta_{n,\varepsilon,X,\theta} > 0$ and let $Y$ be a realization of $\mathbf{Y}$. If there is $l \in (n-1, n]$ such that $W_l(\theta, X, Y) \leq L$, then there is $\tau \in [0, L]$ such that

$$\frac{1}{l} \sum_{y \in Y_\tau^{\tau+l}} f(d(y - \tau, X_0^l)) \geq \theta,$$

which implies that for some $t = k\varepsilon$, $k = 0, 1, \ldots, \lfloor L/\varepsilon \rfloor$,

$$\frac{1}{l} \sup_{\tau \in [t, t+\varepsilon]} \sum_{y \in Y_\tau^{\tau+l}} f(d(y - \tau, X_0^l)) \geq \theta.$$

Since for any $\tau \in [t, t+\varepsilon]$, $Y_{t+\varepsilon}^{t+n-1} \subset Y_\tau^{\tau+l} \subset Y_t^{t+n+\varepsilon}$, the above equality leads to

$$\frac{1}{n-1} \sup_{\tau \in [t, t+\varepsilon]} \sum_{y \in Y_t^{t+n+\varepsilon}} f^+(d(y - \tau, X_0^l))$$

$$- \frac{1}{n} \inf_{\tau \in [t, t+\varepsilon]} \sum_{y \in Y_{t+\varepsilon}^{t+n-1}} f^-(d(y - \tau, X_0^l)) \geq \theta.$$

Because $|d(y - \tau, X_0^l) - d(y - t, X_0^l)| \leq \varepsilon$ for any $y \in \mathbb{R}$ and $\tau \in [t, t+\varepsilon]$, by (2.1), the above inequality implies

$$\frac{1}{n-1} \sum_{y \in Y_t^{t+n+\varepsilon}} f_\varepsilon^+(d(y - t, X_0^l)) - \frac{1}{n} \sum_{y \in Y_{t+\varepsilon}^{t+n-1}} f_\varepsilon^-(d(y - t, X_0^l)) \geq \theta$$

$$\implies B_{n,\varepsilon}(X, Y - t) \geq \theta.$$

Because $t = k\varepsilon$, for some $k = 0, 1, \ldots, \lfloor L/\varepsilon \rfloor$, by the stationarity of $\mathbf{Y}$,

$$\Pr\{\exists l \in (n-1, n] \text{ s.t. } W_l(\theta, X, \mathbf{Y}) \leq L\}$$

$$\leq \Pr\left\{\bigcup_{k=0}^{\lfloor L/\varepsilon \rfloor} \{B_{n,\varepsilon}(X, \mathbf{Y} - k\varepsilon) \geq \theta\}\right\}$$

$$\leq \sum_{k=0}^{\lfloor L/\varepsilon \rfloor} \Pr\{B_{n,\varepsilon}(X, \mathbf{Y} - k\varepsilon) \geq \theta\} = (L/\varepsilon + 1)\beta_{n,\varepsilon,X,\theta}.$$

For $L \geq 1$, this implies

$$\Pr\{\exists l \in (n-1, n] \text{ s.t. } \max(W_l(\theta, X, \mathbf{Y}), 1) \leq L\} \leq 2L\beta_{n,\varepsilon,X,\theta}/\varepsilon.$$

The above bound holds for $L \in (0, 1)$ as well. Choose $L = L(n) = e^{-c(n)}/\beta_{n,\varepsilon,X,\theta}$ with $\sum e^{-c(n)} < \infty$ and $\frac{c(n)}{n} \to 0$ to get

$$\Pr\left\{\exists l \in (n-1, n] \text{ s.t. } \frac{1}{l}\log[\beta_{n,\varepsilon,X,\theta} \times \max(W_l(\theta, X, \mathbf{Y}), 1)] \leq -\frac{c(n)}{l}\right\}$$

$$\leq 2e^{-c(n)}/\varepsilon.$$

By an argument similar to the end of the proof of Lemma 2, (2.3) is proved. $\square$

PROOF OF LEMMA 4. The proof is an application of the Gärtner–Ellis theorem. We will only consider the LDP of $A_n(X, \mathbf{Y})$. The LDP of $B_{n,\varepsilon}(X, \mathbf{Y})$ can be similarly treated.

The first step is to show that almost surely, for $X \sim \mathbf{X}$,

(3.1) $$\frac{1}{n}\log E[e^{ntA_n(X,\mathbf{Y})}] \to \Lambda(t) \quad \text{for all } t \in \mathbb{R}.$$

Let $g_n(y) = \inf_{n-1 \leq l \leq n} f^+(d(y, X_0^l))$ and $h_n(y) = \sup_{n-1 \leq l \leq n} f^-(d(y, X_0^l))$. Then given $t \in \mathbb{R}$,

$$\frac{1}{n}\log E[e^{ntA_n(X,\mathbf{Y})}] = \frac{1}{n}\log E\left[\exp\left\{t\left(\sum_{y \in \mathbf{Y}_0^{n-1}} g_n(y) - \frac{n}{n-1}\sum_{y \in \mathbf{Y}_0^n} h_n(y)\right)\right\}\right]$$

$$= I_1 + I_2,$$

with

$$I_1 = \frac{\lambda}{n}\int_0^{n-1}\left[\exp\left\{t\left(g_n(y) - \frac{n}{n-1}h_n(y)\right)\right\} - 1\right]dy,$$

$$I_2 = \frac{\lambda}{n}\int_{n-1}^n\left[\exp\left\{-\frac{tn}{n-1}h_n(y)\right\} - 1\right]dy.$$

Because $f$ is bounded, $I_2 \to 0$ as $n \to \infty$. Letting $s_n = \min(X_0^{n-1})$ and $\tau_n = \max(X_0^{n-1})$, it is seen that if $s_n \leq y \leq \tau_n$, then $d(y, X_0^l) = d(y, X)$, yielding $g_n(y) = f^+(d(y, X))$ and $h_n(y) = f^-(d(y, X))$. Let

$$F(y) = \exp\left\{t\left(f^+(d(y,X)) - \frac{n}{n-1}f^-(d(y,X))\right)\right\} - 1.$$

Clearly $s_n/n \to 0$. By Proposition 2, we can assume $(n - \tau_n)/n \to 0$. Let $J_n = [0, s_n] \cup [\tau_n, n-1]$. Then by the boundedness of $f$, as $n \to \infty$,

$$I_1 = \frac{\lambda}{n}\int_0^{n-1} F - \frac{\lambda}{n}\int_{J_n} F + \frac{\lambda}{n}\int_{J_n}\left[\exp\left\{t\left(g_n(y) - \frac{n}{n-1}h_n(y)\right)\right\} - 1\right]dy$$

$$= \frac{\lambda}{n}\int_0^n\left[\exp\left\{t\left(f^+(d(0, X-y)) - \frac{n}{n-1}f^-(d(0, X-y))\right)\right\} - 1\right]dy + o(1).$$

Because $\mathbf{X}$ is ergodic, it is seen that $I_1 \to \lambda E[e^{tf(d(0,\mathbf{X}))} - 1]$, proving (3.1) for fixed $t$. It follows that almost surely, (3.1) holds for $t$ in a countable dense subset of $\mathbb{R}$. On the other hand, by the boundedness of $f$, it is not hard to show that $\frac{1}{n}\log E[e^{ntA_n(X,\mathbf{Y})}]$, $n \geq 1$, are equicontinuous functions in $t$ on any bounded region and $\Lambda(t)$ is continuous. Therefore, almost surely, for $X \sim \mathbf{X}$, the convergence in (3.1) holds for all $t \in \mathbb{R}$.

The function $\Lambda(t)$ is smooth and strictly convex. By condition 3 of Theorem 1, $\Lambda(t) \to \infty$ exponentially fast as $t \to \infty$. To finish the proof, consider the event $E = \{f(d(0,\mathbf{X})) < 0\}$. If $\Pr(E) > 0$, then, as $t \to -\infty$, $\Lambda(t) \to \infty$ exponentially fast and hence $\Lambda$ is essentially smooth (cf. [[9]], Definition 2.3.5). By the Gärtner–Ellis theorem, the LDP holds for $A_n(X, \mathbf{Y})$ with the good rate function $\Lambda^*$. If $\Pr(E) = 0$, or equivalently, $f(d(0,\mathbf{X})) \geq 0$ w.p.1, then by Theorem 2.3.6 and Lemma 2.3.9 of [9], for any open set $G$,

$$\liminf_{n \to \infty} \frac{1}{n}\log \Pr\{A_n(X, \mathbf{Y}) \in G\} \geq - \inf_{\alpha \in G \cap (0, \infty)} \Lambda^*(\alpha).$$

Since for $\alpha < 0$, $\Lambda^*(\alpha) = \infty$, and for $0 \leq \alpha < \phi$, $\Lambda^*(\alpha) < \infty$ is decreasing, the above inequality implies

$$\liminf_{n \to \infty} \frac{1}{n}\log \Pr\{A_n(X, \mathbf{Y}) \in G\} \geq - \inf_{\alpha \in G} \Lambda^*(\alpha).$$

Therefore the LDP is proved. $\square$

**4. Waiting times for vector-valued matching scores.** Let comparison or maximization of vectors be made component-wise, for example, if $f = (f_1, \ldots, f_n)$, then $f^+ = (f_1^+, \ldots, f_n^+)$, $\sup_{x \in A} f(x) = (\sup_{x \in A} f_1(x), \ldots, \sup_{x \in A} f_n(x))$, and for vectors $u = (u_1, \ldots, u_n)$, $v = (v_1, \ldots, v_n)$, $\max(u, v) = (\max(u_1, v_1), \ldots, \max(u_n, v_n))$. Given $\theta \in \mathbb{R}^n$, define $W_l(\theta, X, Y)$ as in the case where $f$ is scalar-valued.

PROOF OF THEOREM 2. Lemmas 1–3 still hold. Following the proof for Lemma 4,

$$\frac{1}{n} \log E[e^{n\langle t, A_n(X,\mathbf{Y})\rangle}] \to \Lambda(t). \tag{4.1}$$

Let $\zeta = f(d(0,\mathbf{X}))$. Since $\Lambda(t) < \infty$ on $\mathbb{R}^n$ and is differentiable, to show that the laws of $A_n(X,\mathbf{Y})$ follow the LDP with the good rate function $\Lambda^*(z)$, by the Gärtner–Ellis theorem, it is enough to show that $|\nabla \Lambda(t)| = |E[\zeta e^{\langle t, \zeta \rangle}]| \to \infty$ as $|t| \to \infty$. Assume for a sequence $t_j \in \mathbb{R}^n$ with $|t_j| \to \infty$, $|E[\zeta e^{\langle t_j, \zeta \rangle}]| \le M$. Then there is a subsequence of $\tau_j := t_j/|t_j|$ converging to some $v$ with $|v| = 1$. Without loss of generality, assume the whole sequence $\tau_j$ converges to $v$. Then $|E[\langle v, \zeta \rangle e^{\langle t_j, \zeta \rangle}]| \le M$. By condition $3'$, there is $\varepsilon > 0$ such that $\Pr\{\langle v, \zeta \rangle > 3\varepsilon\} > 0$. Because $f$ is bounded, for $j$ large enough, $|\tau_j - v||\zeta| < \varepsilon$. Then

$$\begin{aligned}
|E[\langle v, \zeta \rangle e^{\langle t_j, \zeta \rangle}]| &\ge E[\langle v, \zeta \rangle e^{\langle t_j, \zeta \rangle} \mathbf{1}_{\{\langle v, \zeta \rangle \ge 3\varepsilon\}}] + E[\langle v, \zeta \rangle e^{\langle t_j, \zeta \rangle} \mathbf{1}_{\{\langle v, \zeta \rangle \le 0\}}] \\
&\ge E[\langle v, \zeta \rangle e^{|t_j|(\langle v, \zeta \rangle - \varepsilon)} \mathbf{1}_{\{\langle v, \zeta \rangle \ge 3\varepsilon\}}] \\
&\quad + E[\langle v, \zeta \rangle e^{|t_j|(\langle v, \zeta \rangle + \varepsilon)} \mathbf{1}_{\{\langle v, \zeta \rangle \le 0\}}] \\
&\ge 3\varepsilon \Pr\{\langle v, \zeta \rangle \ge 3\varepsilon\} e^{2\varepsilon |t_j|} - E|\zeta| e^{\varepsilon |t_j|} \to \infty,
\end{aligned}$$

which is a contradiction.

Let $M(t) = E[e^{\langle t, \zeta \rangle}]$. For any $a > 1$, let $V = \{t : M(t) \le a\}$. Because $M(t)$ is convex and continuous, $V$ is a convex closed set. Assume $V$ is unbounded, then there are $t_j \in V$ with $|t_j| \to \infty$ and $\tau_j = t_j/|t_j| \to v$ for some $v$ with length 1. Given $r > 0$, $|t_j| > r$ for all large $j$. As $0, |t_j|\tau_j \in V$, $r\tau_j \in V$, implying $rv \in V$. As a result, $M(rv) \le a$ for all $r > 0$, which is impossible due to condition $3'$. Therefore, $V$ is bounded. Suppose $|v| \le R$ for all $v \in V$. Then for $t$ with $|t| > R$, by the Hölder inequality, $M(t) \ge (M(Rt/|t|))^{|t|/R} \ge a^{|t|/R}$, and hence $\Lambda(t) = M(t) - 1 \to \infty$ exponentially fast in $|t|$. Therefore, $\Lambda^*(z) \le \sup_{t \in \mathbb{R}}\{|z||t| - \Lambda(t)\}$ is bounded on any bounded set. Since $\Lambda^*$ is convex, then it is seen $\Lambda^*$ is continuous.

By (2.2) and the LDP for the conditional laws of $A_n$, almost surely, for $X \sim \mathbf{X}$ and $Y \sim \mathbf{Y}$,

$$\begin{aligned}
\limsup_{l \to \infty} \frac{1}{l} \log W_l(\theta, X, Y) &\le -\liminf_{n \to \infty} \frac{1}{n} \Pr\{A_n(X, \mathbf{Y}) > \theta\} \\
&\le \inf_{z > \theta} \Lambda^*(z) = \inf_{z \ge \theta} \Lambda^*(z),
\end{aligned} \tag{4.2}$$

with the last equality due to the continuity of $\Lambda^*$. For $z \ge \theta > \phi$, $\langle \mathbf{1}, z \rangle \ge \langle \mathbf{1}, \theta \rangle > \langle \mathbf{1}, \phi \rangle = \lambda E[\langle \mathbf{1}, \zeta \rangle]$. Then by Theorem 1

$$\Lambda^*(z) \ge \sup_{t \ge 0}\{t\langle \mathbf{1}, z \rangle - \lambda E[e^{t\langle \mathbf{1}, \zeta \rangle} - 1]\} \ge \sup_{t \ge 0}\{t\langle \mathbf{1}, \theta \rangle - \lambda E[e^{t\langle \mathbf{1}, \zeta \rangle} - 1]\} > 0.$$

On the other hand, by the Gärtner–Ellis theorem,

$$\liminf_{l\to\infty} \frac{1}{l}\log\max\{1, W_l(\theta, X, Y)\} \geq -\limsup_{n\to\infty} \frac{1}{n}\Pr\{B_{n,\varepsilon}(X, \mathbf{Y}) \geq \theta\}$$
(4.3)
$$\geq \inf_{z\geq\theta} \Lambda_\varepsilon^*(z).$$

Similarly to the proof for Theorem 1, it just remains to show

(4.4) $$\lim_{\varepsilon\to 0}\inf_{z\geq\theta} \Lambda_\varepsilon^*(z) = \inf_{z\geq\theta} \Lambda^*(z),$$

where $\Lambda_\varepsilon^*(z) = \sup_{t\in\mathbb{R}}\{\langle z, t\rangle - E[e^{\langle t, f_\varepsilon(d(0,\mathbf{X}))\rangle} - 1]\}$.

Let $M = \sup_{x\geq 0}|f(x)|$. Since

$$\Lambda_\varepsilon^*(z) \geq \sup_{t\geq 0}\{|z|t - e^{tM}\} \to \infty, \qquad |z|\to\infty,$$

uniformly for $\varepsilon > 0$, for some bounded closed set $A \subset \{z : z \geq \theta\}$, $\inf_{z\geq\theta}\Lambda_\varepsilon^*(z) = \inf_{z\in A}\Lambda_\varepsilon^*(z)$. Next show that as a family of functions parameterized by $\varepsilon > 0$, $\Lambda_\varepsilon^*$ is equicontinuous on $A$ for all small $\varepsilon$.

By the boundedness of $f$ and conditions 2 and 3', for any $v \in \Gamma = \{z : |z| = 1\}$, there are $\eta = \eta(v) > 0$, $\delta = \delta(v) > 0$, and an open neighborhood $U = U(v) \subset \Gamma$, such that

$$\Pr\{\langle v, f_\varepsilon(d(0, \mathbf{X}))\rangle \geq 2\eta\} > \eta \qquad \text{for all } \varepsilon \leq \delta$$

and $M|v - u| < \eta$, for all $u \in U$. Because $\Gamma$ is compact, there are $v_1, \ldots, v_n$ such that $\Gamma = \bigcup_{k=1}^n U(v_k)$. Let $\delta = \min_{k=1}^n \delta(v_k)$ and $\eta = \min_{k=1}^n \eta(v_k)$. For any $v \in \Gamma$, there is $k$ such that $v \in U(v_k)$. Then for any $\varepsilon \leq \delta$, when $\langle v_k, f_\varepsilon(d(0, X))\rangle \geq 2\eta(v_k)$,

$$\langle v, f_\varepsilon(d(0,X))\rangle \geq \langle v_k f_\varepsilon(d(0,X))\rangle - |v - v_k|M \geq 2\eta(v_k) - \eta(v_k) \geq \eta,$$

implying $\Pr\langle v, f_\varepsilon(d(0, \mathbf{X}))\rangle > \eta > \eta$. Fix $L > 0$ such that $|z| \leq L$ for all $z \in A$. For $t \in \mathbb{R} \setminus \{0\}$, write $t = |t|v$. Then as $|t| \to \infty$,

$$\langle z, t\rangle - \lambda E[e^{\langle t, f_\varepsilon(d(0,\mathbf{X}))\rangle} - 1]$$
$$\leq L|t| - \lambda E[(e^{|t|\langle v, f_\varepsilon(d(0,\mathbf{X}))\rangle} - 1)\mathbf{1}_{\{\langle v, f_\varepsilon(d(0,\mathbf{X}))\rangle > \eta\}}] + \lambda$$
$$\leq L|t| - \eta\lambda e^{\eta|t|} + \lambda \to -\infty,$$

uniformly for $z \in A$ and $\varepsilon \leq \delta$. Since $\Lambda_\varepsilon^*(z) \geq 0$, this implies that there is $R > 0$ such that for all $z \in A$ and $\varepsilon \leq \delta$, the maximizer $t^*(z, \varepsilon)$ of $\langle z, t\rangle - \lambda E[e^{\langle t, f_\varepsilon(d(0,\mathbf{X}))\rangle} - 1]$ is in $B_R := \{z : |z| \leq R\}$. Then for any $z_1, z_2 \in A$, it is seen $\Lambda_\varepsilon^*(z_1) - \Lambda_\varepsilon^*(z_2) \leq \langle t^*(z_1, \varepsilon), z_1 - z_2\rangle \leq R|z_1 - z_2|$. Likewise, $\Lambda_\varepsilon^*(z_2) - \Lambda_\varepsilon^*(z_1) \leq R|z_1 - z_2|$. So $\Lambda_\varepsilon^*(z)$ is equicontinuous.

Choose $\varepsilon_n$ such that $\lim_n \inf_{z\geq\theta}\Lambda_n^*(z) = \liminf_{\varepsilon\to 0}\inf_{z\geq\theta}\Lambda_\varepsilon^*(z)$, where $\Lambda_n^* := \Lambda_{\varepsilon_n}^*$. Let $z_n \in A$ be such that $\Lambda_n^*(z_n) = \inf_{z\in A}\Lambda_n^*(z)$. Then $z_n$ has a convergent subsequence. Without loss of generality, suppose $z_n \to z \in A$. Following the same argument as in the proof of Theorem 1, $\Lambda_n^*(z) \to \Lambda^*(z)$. Then by the equicontinuity of $\Lambda_\varepsilon^*$,

$$\liminf_{\varepsilon \to 0} \inf_{z \geq \theta} \Lambda_\varepsilon^*(z) = \lim_{n \to \infty} \Lambda_n^*(z_n) = \lim_{n \to \infty} \Lambda_n^*(z) = \Lambda^*(z) \geq \inf_{z \geq \theta} \Lambda^*(z) > 0.$$

Therefore, (4.3) can be replaced by

$$\liminf_{l \to \infty} \frac{1}{l} \log W_l(\theta, X, Y) \geq \inf_{z \geq \theta} \Lambda_\varepsilon^*(z).$$

This together with (4.2) implies that

$$\limsup_{\varepsilon \to 0} \inf_{z \geq \theta} \Lambda_\varepsilon^*(z) \leq \inf_{z \geq \theta} \Lambda^*(z),$$

which completes the proof of (4.4). □

**5. An approximation for large deviations.** Given $\theta > \phi := \lambda E[f(d(0, \mathbf{X}))]$, it is easy to see that $\theta t - \Lambda(t)$ achieves $\Lambda^*(\theta)$ at a unique point $t_0$. Furthermore, $t_0 \in (0, \infty)$ and

$$(5.1) \qquad \theta = \Lambda'(t_0) = \lambda E[f(d(0, \mathbf{X}))e^{t_0 f(d(0, \mathbf{X}))}].$$

LEMMA 5. *Almost surely, for $X \sim \mathbf{X}$, when $l$ is large, $\theta t - \Lambda_{X_0^l, l}(t)$ achieves $\Lambda_{X_0^l, l}^*(\theta)$ on $(0, \infty)$ and the maximizer $t^* = t^*(X, l)$ is unique. Furthermore, $t^*$ satisfies*

$$(5.2) \qquad \theta = \Lambda'_{X_0^l, l}(t^*) = \frac{\lambda}{l} \int_0^l f(d(y, X_0^l)) e^{t^* f(d(y, X_0^l))} \, dy$$

*and, as $l \to \infty$, $t^* \to t_0$, $\Lambda_{X_0^l, l}(t^*) \to \Lambda(t_0)$ and $\Lambda''_{X_0^l, l}(t^*) \to \Lambda''(t_0)$.*

PROOF. Almost surely, for $X \sim \mathbf{X}$, for all large $l$, $\Lambda_{X_0^l, l}(t)$ is smooth, strictly convex, $\Lambda_{X_0^l, l}(0) = 0$, and $\Lambda'_{X_0^l, l}(t) \to \infty$ exponentially fast as $t \to \infty$. Furthermore, following the proof of Lemma 1,

$$\lim_{l \to \infty} \frac{\lambda}{l} \int_0^l f(d(y, X_0^l)) \, dy = \lim_{l \to \infty} \frac{\lambda}{l} \int_0^l f(d(y, X)) \, dy = E[f(d(0, \mathbf{X}))].$$

and hence $\Lambda'_{X_0^l, l}(0) < \theta$, implying $\theta t - \Lambda_{X_0^l, l}(t)$ has a unique maximizer $t^*$ which is in $(0, \infty)$. By differentiation, (5.2) is proved. For any $t > t_0$, by (5.1) and (5.2), as $l \to \infty$,

$$\Lambda'_{X_0^l, l}(t) \to \Lambda'(t) > \Lambda'(t_0) = \theta = \Lambda'_{X_0^l, l}(t^*).$$

Therefore, $t^* < t$ eventually, giving $\limsup_{l \to \infty} t^* \leq t_0$. Likewise, $\liminf_{l \to \infty} t^* \geq t_0$. This proves $t^* \to t_0$. Finally, following the equicontinuity argument as in the previous sections,

$$\Lambda_{X_0^l, l}(t^*) = \frac{\lambda}{l} \int_0^l [e^{t^* f(d(y, X_0^l))} - 1] \, dy \to \lambda E[e^{t_0 f(d(0, \mathbf{X}))} - 1] = \Lambda(t_0)$$

and

$$\Lambda''_{X_0^l, l}(t^*) = \frac{\lambda}{l} \int_0^l f^2(d(y, X_0^l)) e^{t^* f(d(y, X_0^l))} \, dy$$

$$\to \lambda E[f^2(d(0, \mathbf{X})) e^{t_0 f(d(0, \mathbf{X}))}] = \Lambda''(t_0) > 0. \qquad \square$$

PROOF OF THEOREM 3. Given $X \sim \mathbf{X}$ such that $\Lambda_{X_0^l, l}$ has the properties described in Lemma 5, let

$$J_l = \exp(l \Lambda^*_{X_0^l, l}(\theta)) \Pr \left\{ \sum_{y \in \mathbf{Y}_0^l} f(d(y, X_0^l)) \geq l\theta \right\}.$$

First, because $J_l \leq \exp(l \Lambda^*_{X_0^l, l}(\theta)) E[\exp\{t^*(\sum_{y \in \mathbf{Y}_0^l} f(d(y, X_0^l)) - l\theta)\}] = 1$, we have

$$\limsup_{l \to \infty} \frac{1}{\sqrt{l}} \log J_l \leq 0.$$

It remains to show that

(5.3) $$\liminf_{l \to \infty} \frac{1}{\sqrt{l}} \log J_l \geq 0.$$

For $l > 0$ large enough, let $g(y) := f(d(y, X_0^l))$. Let $t^* > 0$ be the maximizer of $\theta t - \Lambda_{X_0^l, l}(t)$ as in Lemma 5. Define measures $\nu = \nu_{X_0^l, l}$ and $\mu = \mu_{X_0^l, l}$ on $[0, l]$, respectively, by

$$\frac{d\nu(y)}{dy} = \lambda e^{t^* g(y)} \quad \text{and} \quad d\mu(y) = \frac{d\nu(y)}{K},$$

with $K = \int_0^l d\nu(y)$. Then $\mu$ is a probability measure. It is easy to see that

(5.4) $\quad K = l(\Lambda_{X_0^l, l}(t^*) + \lambda) = l(\theta t^* - \Lambda^*_{X_0^l, l}(\theta) + \lambda) \quad \text{and} \quad l\theta = K E[g(\xi)],$

with $\xi \sim \mu$. Also,

(5.5) $\quad \dfrac{K}{l} = \Lambda_{X_0^l, l}(t^*) + \lambda \to \Lambda(t_0) + \lambda = \lambda E[e^{t_0 f(d(0, \mathbf{X}))}] > 0 \qquad \text{as } l \to \infty.$

Letting $m = E[g(\xi)]$, by (5.4) and the properties of Poisson processes,

$$J_l = e^{l\Lambda^*_{X_0^l,l}(\theta)} \sum_{n=0}^{\infty} e^{-\lambda l} \frac{\lambda^n}{n!} \int_{[0,l]^n} \mathbf{1}_{\{\sum_{i=1}^n g(y_i) \geq l\theta\}} dy_1 \cdots dy_n$$

$$= e^{l\Lambda^*_{X_0^l,l}(\theta) - \lambda l + K} \sum_{n=0}^{\infty} e^{-K} \frac{K^n}{n!} \int_{[0,l]^n} \mathbf{1}_{\{\sum_{i=1}^n g(y_i) \geq l\theta\}} \exp\left\{-t^* \sum_{i=1}^n g(y_i)\right\}$$

$$\times \prod_{i=1}^n \frac{\lambda e^{t^* g(y_i)}}{K} dy_1 \cdots dy_n$$

$$= e^{l\theta t^*} \sum_{n=0}^{\infty} e^{-K} \frac{K^n}{n!} E\left[\mathbf{1}_{\{\sum_{i=1}^n g(\xi_i) \geq l\theta\}} \exp\left\{-t^* \sum_{i=1}^n g(\xi_i)\right\}\right]$$

$$= \sum_{n=0}^{\infty} e^{-K} \frac{K^n}{n!} E\left[\mathbf{1}_{\{\sum_{i=1}^n (g(\xi_i) - m) \geq (K-n)m\}} \exp\left\{-t^* \sum_{i=1}^n (g(\xi_i) - m)\right\}\right] e^{t^*(K-n)m},$$

with $\xi_1, \ldots, \xi_n$ i.i.d. $\sim \mu$.

Fix $\delta > 0$. Recall $t^* > 0$. If $m \geq 0$, then

$$J_l \geq \sum_{K \leq n \leq K + \delta\sqrt{K}} e^{-K} \frac{K^n}{n!} E\left[\mathbf{1}_{\{\sum_{i=1}^n (g(\xi_i) - m) \geq 0\}} \exp\left\{-t^* \sum_{i=1}^n (g(\xi_i) - m)\right\}\right]$$

$$\times e^{-t^* \sqrt{K} m \delta}.$$

A similar bound can be obtained when $m < 0$, by summing over $K - \delta\sqrt{K} \leq n \leq K$ instead. Without loss of generality, assume $m \geq 0$. Let

$$G_n = \frac{\sum_{i=1}^n (g(\xi_i) - m)}{\sqrt{n \operatorname{Var}[g(\xi)]}}.$$

Let $l \to \infty$. Then $t^* \to t_0$ and by (5.5), $K \to \infty$. There is a constant $c_1 = c_1(\delta) > 0$, such that for large $K$, $\sum_{K \leq n \leq K + \delta\sqrt{K}} e^{-K} K^n/n! \geq c_1$ and hence

$$J_l \geq \sum_{K \leq n \leq K + \delta\sqrt{K}} e^{-K} \frac{K^n}{n!} E[\mathbf{1}_{\{0 \leq G_n \leq \delta\}} e^{-t^* \sqrt{n \operatorname{Var}[g(\xi)]} G_n}] e^{-t^* \sqrt{K} m \delta}$$

$$\geq \sum_{K \leq n \leq K + \delta\sqrt{K}} e^{-K} \frac{K^n}{n!} \Pr\{0 \leq G_n \leq \delta\} e^{-t^* \sqrt{2K \operatorname{Var}[g(\xi)]} \delta} e^{-t^* \sqrt{K} m \delta}$$

$$\geq c_1 \min_{K \leq n \leq K + \delta\sqrt{K}} \Pr\{0 \leq G_n \leq \delta\} e^{-t^* \sqrt{K} D \delta}$$

with $D = \sqrt{2 \operatorname{Var}[g(\xi)]} + m$. It is not hard to see that for $\xi \sim \mu$ and $\eta = f(d(0, \mathbf{X}))$,

$$\operatorname{Var}[g(\xi)] = \frac{\int_0^l g^2(y) e^{t^* g(y)} dy}{\int_0^l e^{t^* g(y)} dy} - \left(\frac{\int_0^l g(y) e^{t^* g(y)} dy}{\int_0^l e^{t^* g(y)} dy}\right)^2$$

$$\to \frac{E[\eta^2 e^{t_0 \eta}]}{E[e^{t_0 \eta}]} - \left(\frac{E[\eta e^{t_0 \eta}]}{E[e^{t_0 \eta}]}\right)^2 > 0$$

and hence $\mathrm{Var}[g(\xi)]$ is uniformly bounded below from 0 for all large $l$. Because $g(y) = f(d(y, X_0^l))$ is uniformly bounded, $G_n$ satisfy Lindeberg's condition, giving $G_n \xrightarrow{D} N(0,1)$. Together with (5.5), these imply that there is a constant $c_2 > 0$ which is independent of $l$ and $\delta$, and some $\rho = \rho(\delta) > 0$, such that $J_l \geq \rho e^{-c_2 t^* \sqrt{l} \delta}$, yielding

$$\liminf_{l \to \infty} \frac{1}{\sqrt{l}} \log J_l \geq -c_2 t_0 \delta.$$

Because $\delta$ is arbitrary, (5.3) is proved. $\square$

## 6. Asymptotic normality.

PROOF OF THEOREM 4. From Lemma 5, it is seen that almost surely, for large $l > 0$, there are unique $\tau_l, t_l > 0$ with $\Lambda^*_{X,l} = \theta \tau_l - \Lambda(\tau_l)$, $\Lambda^*_{X_0^l, l}(\theta) = \theta t_l - \Lambda_{X_0^l, l}(t_l)$. Furthermore, $\tau_l, t_l \to t_0$ as $l \to \infty$. Fix $\delta, M > 0$, such that $\tau_l, t_l \in (t_0 - \delta, t_0 + \delta)$ for all large $l > 0$ and $\lambda |e^{tf(d(y, X_0^l))} - 1| \leq M/2$ for $t \in (t_0 - \delta, t_0 + \delta)$ and all $y$. Following the argument in the proof of Lemma 1, on $(t_0 - \delta, t_0 + \delta)$,

$$(6.1) \qquad |\Lambda_{X,l}(t) - \Lambda_{X_0^l, l}(t)| \leq (\min(X_0^l) + d_l) M/l,$$

where $d_l = l - \max(X_0^l)$. Clearly, $\min(X_0^l) = O(1)$ w.p.1. Letting $n = \lfloor l \rfloor$, $d_l \leq s_n = n + 1 - \max(X_{-\infty}^n) \stackrel{D}{=} 1 - \max(X_{-\infty}^0) \stackrel{D}{=} 1 + \rho U$, with $U \sim \mathrm{Exp}(1)$. Given $\varepsilon > 0$, $\Pr\{s_n \geq \sqrt{\varepsilon n}\} \leq \Pr\{(U+1)^2 \geq \varepsilon n\}$. Since $EU^2 < \infty$, applying the Borel–Cantelli lemma to $s_n$, it is seen that $d_l = o(\sqrt{l})$, w.p.1, and hence the left-hand side of (6.1) is $o(1/\sqrt{l})$ w.p.1. Then, by $t_l, \tau_l \in (t_0 - \delta, t_0 + \delta)$,

$$(6.2) \quad |\Lambda^*_{X,l}(\theta) - \Lambda^*_{X_0^l, l}(\theta)| \leq \sup_{|t - t_0| < \delta} |\Lambda_{X,l}(t) - \Lambda_{X_0^l, l}(t)| = o(1/\sqrt{l})$$

w.p.1.

On the other hand, it is easy to see that for large $l > 0$ and $n = \lfloor l \rfloor$,

$$(6.3) \qquad |\Lambda_{X,l}(t) - \Lambda_{X,n}(t)| \leq 2M/l$$

for all $t \in (t_0 - \delta, t_0 + \delta)$. In particular, letting $t = \tau_n, \tau_l$ leads to

$$(6.4) \quad |\Lambda^*_{X,l}(\theta) - \Lambda^*_{X,n}(\theta)| \leq \sup_{|t - t_0| < \delta} |\Lambda_{X_0^l, l}(t) - \Lambda_{X_0^l, l}(t)| \leq 2M/l.$$

From (6.1)–(6.4), it is seen that (1.8) holds if we can show

$$z_n = \sqrt{n}\,[\theta(\tau_n - t_0) - (\Lambda_{X,n}(\tau_n) - \Lambda_{X,n}(t_0))] = o(1) \qquad \text{w.p.1.}$$

Since $z_n = \Lambda_{X,n}^*(\theta) - (\theta t_0 - \Lambda_{X,n}(t_0)) \geq 0$, it is enough to prove $\limsup z_n \leq 0$, or equivalently,

(6.5) $\qquad \liminf_{n \to \infty} \sqrt{n}\,[\Lambda_{X,n}(\tau_n) - \Lambda_{X,n}(t_0) - \theta(\tau_n - t_0)] \geq 0.$

Because $\tau_n \in (t_0 - \delta, t_0 + \delta)$ and $\Lambda_{X,n}(t)$ is smooth, by Taylor's expansion, for some $t^* \in (t_0 - \delta, t_0 + \delta)$,

$$\Lambda_{X,n}(t) - \Lambda_{X,n}(t_0) - \theta(t - t_0) = A_n(t - t_0) + \frac{B_{n,t^*}(t - t_0)^2}{2} \geq -\frac{A_n^2}{2B_{n,t^*}},$$

where

$$A_n = \frac{1}{n} \int_0^n f(d(y, X)) e^{t_0 f(d(y,X))}\, dy - \theta,$$

$$B_{n,t} = \frac{1}{n} \int_0^n f^2(d(y, X)) e^{t f(d(y,X))}\, dy > 0.$$

Because $f$ is bounded and $\mathbf{X}$ is ergodic, there exists a constant $\eta > 0$, such that $B_{n,t} > \eta$ for all large $n$ and $t \in (t_0 - \delta, t_0 + \delta)$. The random variables

$$Z_n = \int_{n-1}^n f(d(y, \mathbf{X})) e^{t_0 f(d(y,\mathbf{X}))}\, dy$$

are bounded and form a stationary process such that $A_n = \frac{1}{n}\sum_{k=1}^n Z_k - \theta$. Since $t_0$ maximizes $\theta t - E[e^{tf(d(0,\mathbf{X}))}]$, $\theta = E[f(d(0,\mathbf{X}))e^{t_0 f(d(0,\mathbf{X}))}] = EZ_n$. Let $\alpha(k) := \sup\{|P(F_1 \cap F_2) - P(F_1)P(F_2)| : F_1 \in \sigma(Z_n, n \leq m),\ F_2 \in \sigma(Z_n, n > m+k),\ m \geq 1\}$. We shall show $\sum_{k=1}^\infty \alpha(k) < \infty$, once this is done, it follows that $\sqrt{n} A_n^2 \to 0$ almost surely (cf. [12], Theorem 2). Then the left-hand side of (6.5) is bounded below by $\liminf(-\sqrt{n}A_n^2/2\eta) = 0$, which completes the proof of (1.8).

Given $k \geq 1$, for any $m \geq 1$, let

$$I = \mathbf{1}_{\{\mathbf{X} \cap (m, m+k/3) \neq \varnothing\}} \quad \text{and} \quad J = \mathbf{1}_{\{\mathbf{X} \cap (m+2k/3, m+k) \neq \varnothing\}}.$$

From the definition of $Z_n$, it is seen that when $I = 1$, for $n \leq m$, $Z_n$ only depends on $\mathbf{X}_{-\infty}^{m+k/3}$. Therefore, for any event $F_1 \in \sigma(Z_n, n \leq m)$, $F_1 \cap \{I = 1\} \in \sigma(\mathbf{X}_{-\infty}^{m+k/3})$. Likewise, for any event $F_2 \in \sigma(Z_n, n > m+k)$, $F_2 \cap \{J = 1\} \in \sigma(\mathbf{X}_{m+2k/3}^\infty)$. Consequently, by the property of Poisson processes, $P(F_1 \cap F_2, I = 1, J = 1) = P(F_1, I = 1)$. Because $\mathbf{X}$ is stationary and has density $\rho > 0$,

$$0 \leq P(F_1 \cap F_2) - P(F_1 \cap F_2, I = 1, J = 1)$$
$$\leq P\{I = 0\} + P\{J = 0\} = 2P\{\mathbf{X} \cap [0, k/3] = \varnothing\} = 2e^{-\rho k/3}$$

and similarly, $0 \leq P(F_1)P(F_2) - P(F_1, I = 1)P(F_2, J = 1) \leq 2e^{-\rho k/3}$. Therefore, $|P(F_1 \cap F_2) - P(F_1)P(F_2)| \leq 4e^{-\rho k/3}$, leading to $\sum \alpha(k) \leq 4\sum e^{-\rho k/3} < \infty$.

By (6.1) and (6.3), in order to show (1.9), it is enough to demonstrate $\sqrt{n}(\theta t_0 - \Lambda_{X,n}(t_0) - \Lambda^*(\theta)) \xrightarrow{D} N(0, 4\rho\sigma^2)$, or

$$\frac{1}{\sqrt{n}}\left(\int_0^n g(d(y, \mathbf{X}))\,dy - n\nu\right) \xrightarrow{D} N(0, 4\rho\sigma^2),$$

where $\nu = E[g(d(0, \mathbf{X}))]$. Because $d(0, \mathbf{X}) \sim \frac{1}{2\rho}U$, with $U \sim \text{Exp}(1)$,

$$\nu = E\left[g\left(\frac{U}{2\rho}\right)\right] = 2\rho E\left[G\left(\frac{U}{2\rho}\right)\right].$$

For $\mathbf{X}_0^n = \{x_1, \ldots, x_N\}$, with $x_i < x_{i+1}$, letting $x_0 = 0$, $x_{N+1} = n$ and $I = \sum_{i=0}^N G(\frac{x_{i+1}-x_i}{2})$,

$$\int_0^n g(d(y, \mathbf{X}))\,dy$$

$$= \int_0^{x_1} g(y)\,dy + 2\sum_{i=1}^{N-1}\int_{x_i}^{(x_i+x_{i+1})/2} g(y-x_i)\,dy + \int_{x_N}^n g(y-x_N)\,dy$$

$$= 2I + G(x_1) + G(n - x_N) - 2G(x_1/2) - 2G((n-x_N)/2).$$

The last four terms are $o(n^{-1/2})$, so it suffices to consider $2I$. Given a specific value of $N$,

$$(x_1, x_2, \ldots, x_N) \sim \left(\frac{nU_0}{(N+1)\overline{U}}, \frac{n(U_0+U_1)}{(N+1)\overline{U}}, \ldots, \frac{n(U_0+U_1+\cdots+U_N)}{(N+1)\overline{U}}\right),$$

with $U_0, \ldots, U_N$ i.i.d. $\sim \text{Exp}(1)$, and $\overline{U}_N = \frac{1}{N+1}\sum_{k=0}^N U_k$. So by Taylor's expansion,

$$I \stackrel{D}{=} \sum_{i=0}^N G\left(\frac{nU_i}{2(N+1)\overline{U}_N}\right)$$

$$= \sum_{i=0}^N G\left(\frac{U_i}{2\rho}\right) + \sum_{i=0}^N g\left(\frac{n(1-\xi)U_i}{2(N+1)\overline{U}_N} + \frac{\xi U_i}{2(N+1)\overline{U}_N}\right)\left[\frac{n\rho}{(N+1)\overline{U}_N} - 1\right]\frac{U_i}{2\rho}$$

$$= \sum_{i=0}^N G\left(\frac{U_i}{2\rho}\right) + (N+1)A_N\left[\frac{n\rho}{(N+1)\overline{U}_N} - 1\right]$$

$$= \sum_{i=0}^N \left[G\left(\frac{U_i}{2\rho}\right) - \frac{A_N}{\overline{U}_N}(U_i - 1)\right] + \frac{A_N}{\overline{U}_N}(n\rho - N - 1),$$

where $\xi \in (0,1)$ and

$$A_N = \frac{1}{N+1}\sum_{i=0}^N g\left(\frac{n(1-\xi)U_i}{2(N+1)\overline{U}_N} + \frac{\xi U_i}{2(N+1)\overline{U}_N}\right)\frac{U_i}{2\rho}.$$

Therefore,

$$\frac{1}{\sqrt{n}}\left(I - \frac{n\nu}{2}\right) \stackrel{D}{=} \frac{1}{\sqrt{n}}\sum_{i=0}^{N}\left[G\left(\frac{U_i}{2\rho}\right) - \frac{\nu}{2\rho} - \frac{A_N}{\overline{U}_N}(U_i - 1)\right]$$
$$+ \frac{N + 1 - n\rho}{\sqrt{n}}\left(\frac{\nu}{2\rho} - \frac{A_N}{\overline{U}_N}\right).$$

As $n \to \infty$, $(N + 1 - n\rho)/\sqrt{n\rho} \stackrel{D}{\to} N(0,1)$. And as $m \to \infty$, $\overline{U}_m \stackrel{P}{\to} 1$, $A_m \stackrel{P}{\to} E[g(\frac{U}{2\rho})\frac{U}{2\rho}]$ (because $g$ is continuous). These combined with CLT then give (1.10).

$\square$

DEPARTMENT OF STATISTICS
UNIVERSITY OF CHICAGO
5734 UNIVERSITY AVENUE
CHICAGO, ILLINOIS 60637
USA
E-MAIL: chi@galton.uchicago.edu
URL: http://galton.uchicago.edu/~chi